\documentclass[10pt,
final,
 journal, 
twocolumn,
 comsoc,
 letterpaper]{IEEEtran} 

\usepackage{srcltx}
\usepackage[psamsfonts]{amsfonts}
\usepackage{amsmath}
\usepackage{amssymb}
\usepackage{amsfonts}
\usepackage{graphics}
\usepackage{tikz} 
\usetikzlibrary{fit,positioning}
\usetikzlibrary{arrows,positioning,automata}
\usepackage{psfrag}
\usepackage{epsfig}
\usepackage{subcaption}
\usepackage{array}
\usepackage{algorithm}
\usepackage{algorithmic}
\usepackage{pifont}
\usepackage{srcltx}
\usepackage[psamsfonts]{amsfonts}
\usepackage{makeidx}  
\usepackage{bbm}
\usepackage{hhline}
\usepackage{eufrak}
\usepackage{yfonts}
\usepackage{color}
\usepackage{url}
\usepackage{dsfont}
\usepackage{caption}
\usepackage[square, comma, sort&compress, numbers]{natbib}
\usepackage{balance} 
\usepackage{amsthm}
\usepackage{verbatim}

\newtheorem{remark}{Remark}

\newtheorem{define}{Definition}
\usepackage{comment}
\newtheorem{theorem}{Theorem}[section]
\newtheorem{lemma}[theorem]{Lemma}

\begin{document}
\pagenumbering{arabic}
\title{\huge{Dynamic Sensor Selection for Reliable Spectrum Sensing via E-optimal Criterion}}
\author{Mohsen Joneidi, Alireza Zaeemzadeh, and Nazanin Rahnavard\\
Department of Electrical and Computer Engineering, University of Central Florida\\
\{joneidi, zaeemzadeh, nazanin\}@eecs.ucf.edu
\vspace{-2mm}

\thanks{This material is based upon work supported by the National Science Foundation under Grant No. CCF-1718195.}}
\maketitle
\begin{abstract}
Reliable and efficient spectrum sensing through dynamic selection of a subset of spectrum sensors is studied. The problem of selecting $K$ sensor measurements from a set of $M$ potential sensors is considered where $K\ll M$. In addition, $K$ may be less than the dimension of the unknown variables of estimation. Through sensor selection, we reduce the problem to an under-determined system of equations with potentially infinite number of solutions. However, the sparsity of the underlying data facilitates limiting the set of solutions to a unique solution. Sparsity enables employing the emerging compressive sensing technique, where the compressed measurements are selected from a large number of potential sensors. This paper suggests selecting sensors in a way that the reduced system of equations constructs a well-conditioned measurement matrix. Our criterion for sensor selection is based on E-optimalily, which is highly related to the restricted isometry property that provides some guarantees for sparse solution obtained by $\ell_1$ minimization. Moreover, the proposed framework exploits a feedback mechanism to evolve the selected sensors dynamically over time. The evolution aims to maximize the reliability of the sensed spectrum.  
\end{abstract}


\section{Introduction}
\vspace{-1mm}
In the last decade, complex systems containing very large numbers of data-gathering devices were developed. An example is wireless sensor networks. In such systems, the processing unit has to deal with an excessively large number of observations acquired by the various sensors. Often there exist some redundancies within the sensed data and they should be pruned. Sensor selection and sensor scheduling aim to address this problem. In many applications the sensor selection task is non-trivial and possibly consists of addressing an NP-hard problem (i.e., there are $M\choose K$ possibilities of choosing $K$ distinct sensors out of $M$ available ones). This essentially implies that an optimal solution cannot be efficiently computed, in particular when the number of sensors becomes excessively large. A convex relaxation of the original NP-hard problem has been suggested in \cite{joshi_ss}. The most prominent advantage of this approach over other methods is its practicality, thanks to many well-established computationally-efficient convex optimization techniques. In addition to convex relaxation, a sub-modular cost function as the criterion of sensor selection allows us to take advantage of greedy optimization methods for selecting sensors \cite{sub_modularity,shirazi}.

The existing studies on sensor selection mostly consider heuristic approaches . For example, in \cite{joshi_ss} the volume of the reduced bases is considered. This method is called \emph{D-optimality}. In addition, \emph{A-optimality} \cite{Boyd:2004} and \emph{E-optimality} \cite{Boyd:2004} are suggested  as some other alternative heuristics already introduced in convex optimization. These heuristics are presented   without any specific justification for sensor selection application. In this paper we are going to exploit a criteria more judiciously in favor of compressed sensing (CS) theoretical guarantees.

Inspired by the compressed sensing theory, this paper suggests to design and optimize a sensor selection method. The goal is to reduce a measurement matrix to only a small fraction of its rows in order to optimize the proposed E-optimal criterion.

The main contributions of the paper are summarized as follows:
\begin{itemize}
\item The link between matrix subset selection, especially volume sampling and sensor selection is investigated.
\item The E-optimal criterion for matrix subset selection is proposed, which results in a new sensor selection method,
\item The suitability of the E-optimal criterion is discussed which is an upper bound for RIP coefficients for compressive sensing, and
\item The reliability concept for power spectrum map is introduced and  it is exploited for reliable sensor selection.
\end{itemize}

Table \ref{tbl_notations} presents the employed notations throughout this paper.
\begin{table}[b]
\centering
\caption{ \small{Employed notations and variables in this paper.}}
\label{tbl_notations}
    \begin{tabular}{  l | l }
    \hline
       Variable Type & Notation \\ \hline\hline
    Constant Scalar & $X$\\ \hline
    Vector & $\boldsymbol{x}$\\ \hline
    $s^{\text{th}}$ entry of Vector & $x_s$\\ \hline
    Matrix & $\boldsymbol{X}$\\ \hline
    Set & $\mathbb{X}$  \\ \hline
    Selected Rows of  $\boldsymbol{A}$ by set $\mathbb{X}$ & $\boldsymbol{A}_{\mathbb{X}}$  \\ \hline
    Number of non-zero entries of $\boldsymbol{x}$ & $\|\boldsymbol{x}\|_0$  \\ \hline
    Trace of Matrix $\boldsymbol{X}$ & Tr($\boldsymbol{X}$)\\ \hline
        Projection of $\boldsymbol{X}$ on its rows set $\mathbb{T}$& $\pi_{\mathbb{T}}(\boldsymbol{X})$\\ \hline
        M & number of potential sensors\\ \hline
            K & number of selected sensors\\ \hline
    \end{tabular}
\end{table}

The rest of paper is organized as follows. Section \ref{problem_state} states the problem of sensor selection and reviews some existing methods. E-optimal sampling is introduced in Section \ref{eopt_sec} and a new sensor selection method is proposed. Section \ref{rel_section} propose dynamic sensor selection based on the reliability. The optimization method for solving the proposed problem is explained in Section \ref{opt_sec}. Section \ref{Experimental} presents the simulation results  and Section \ref{conc} concludes the paper.

\section{Background}
This paper address a joint framework for spectrum sensing with partial sensing from a big set of sensors. The partial data are selected through a sensor selection procedure. Viewing spectrum sensing and sensor selection together in a joint problem is resulted to inspiring theoretical results in addition to a new application.   

The prerequisite background of the proposed framework is reviewed in this section. The first subsection reviews compressed sensing theory and then the system model of spectrum sensing is introduced in the second subsection. The third subsection introduces the sensor selection problem. The last subsection review theoretical results in matrix subset selection literature which is highly related to our proposed selection method.
\subsection{Compressed Sensing}
Compressed sensing is a technique by which sparse signals can be measured at a rate less than conventional Nyquist sampling theorem. \cite{cs_book}. There exist vast applications of CS in signal and image processing \cite{cs_app_iamge}, channel estimation \cite{cs_app_comm} and spectrum sensing \cite{cs_app_ss}. CS aims to recover a sparse vector, $\mathbf{x}$,  using a small number of measurements $\mathbf{y}$. The CS problem can be formulated as,
\begin{equation}
\hat{\textbf{x}}=\underset{\mathbf{x}}{\text{argmin}} \|\mathbf{x}\|_0 \quad \text{s.t.} \; \; \mathbf{y}=\boldsymbol{\Phi} \mathbf{x},
\end{equation}
where, $\Vert . \Vert_0$ represents the number of non-zero elements of a vector. $\boldsymbol{\Phi}\in \mathbb{R}^{K\times N}$ is called measurement matrix that provides us $K$ measurements collected in $\boldsymbol{y}$. Exact solution of the above optimization problem is through the combinational search among all possible  subsets. Due to its high computational burden, this algorithm is impractical for high dimension scenarios. Many sub-optimal algorithms have been proposed such as OMP \cite{omp}, smoothed $\ell_0$ \cite{slzero} and basis pursuit \cite{basispursut}. Basis pursuit is based on relaxing $\ell_0$ to $\ell_1$ norm and is popular due to  theoretical guarantees and reasonable computational burden \cite{dic_learning}
. The theoretical guarantees for $\ell_1$ minimization arise from several sufficient conditions based on some suggested metrics. These include the mutual coherence \cite{Donoho01uncertainty}, null space property \cite{Cohen09compressedsensing}, spark \cite{donoho2003optimally} and restricted isometery property (RIP) \cite{Candes_RIP}. Except for the mutual coherence, none of these measures can be efficiently calculated for an arbitrary given measurement matrix $\boldsymbol{\Phi}$. For example, the RIP requires enumerating over an exponential number of index sets. RIP is defined as follows.

\begin{define}\cite{Candes_RIP}
A measurement matrix is said to satisfy symmetric form RIP of order $S$  with constant $\delta_S$ if $\delta_S$ is the smallest number that\\
\small{
\begin{equation}\label{RIP_def}
(1-\delta_{S})\Vert \textbf{x}\Vert_2^2\le \Vert \boldsymbol{\Phi}\textbf{x}\Vert_2^2\le (1+\delta_{S})\Vert \textbf{x}\Vert_2^2,
\end{equation}
}
holds for every $S$-sparse $\textbf{\text{x}}$ (i.e. $\textbf{\text{x}}$ contains at most $S$ nonzero entries). 
\end{define}
Based on this definition several guarantees are proposed in terms of $\delta_{2S}$, $\delta_{3S}$ and $\delta_{4S}$ in \cite{Candes_RIP2} and \cite{Davies_RIP3} in order to guarantee recovering $S$-sparse vectors. By $S$-sparse we mean a vector that has $S$ non-zero entries.
In \cite{Blanchard:2011} an asymmetric form of definition 1 is introduced in order to more precisely quantify the RIP.
\begin{define} \cite{Blanchard:2011}
\label{def_rip}
For a measurement matrix the asymmetric RIP constants $\delta^L_{S}$ and $\delta^U_S$ are defined as,
\small{
\begin{equation}
\label{asymm_RIP}
\begin{split}
&\delta^L_{S}(\boldsymbol{\Phi})= \underset{c>0}{\text{argmin}}\;  (1-c)\Vert \textbf{x}\Vert_2^2\le \Vert \boldsymbol{\Phi}\textbf{x}\Vert_2^2,\;\; \forall \textbf{x}\in \mathcal{X}_S^N , \\
&\delta^U_{S} (\boldsymbol{\Phi})= \underset{c>0}{\text{argmin}}\;  (1+c)\Vert \textbf{x}\Vert_2^2\geq \Vert \boldsymbol{\Phi}\textbf{x}\Vert_2^2, \;\; \forall \textbf{x}\in \mathcal{X}_S^N,
\end{split}
\end{equation}
}
\end{define}
where, $\mathcal{X}_S^N$ refers to the set of $S$-sparse vectors in $\mathbb{R}^N$.

\begin{remark}\cite{Blanchard:2011}
Although both the smallest and largest singular values of $\boldsymbol{\Phi_\mathbb{S}}^T \boldsymbol{\Phi_\mathbb{S}}$ \footnote{$\mathbb{S}$ represents a set with cardinality of $S$ and $\boldsymbol{\Phi_\mathbb{S}}$ represents the corresponding selected rows of $\boldsymbol{\Phi}$.} affect the stability of the reconstruction algorithms, the smaller eigenvalue is dominant for compressed sensing in that it allows distinguishing between sparse vectors, $\mathcal{X}_S^N$, given their measurements by $\boldsymbol{\Phi}$.
\end{remark}

\subsection{Spectrum sensing problem statement}\label{spec_sens}
Cognitive radio (CR) is a promising solution to alleviate today's spectrum deficiency caused by an increased demand for the wireless technologies \cite{Akyildiz06_xG, Dynamic_Spectrum_Access4205091}.The CR paradigm allows the unlicensed or secondary users (SUs) to coexist with the PUs. The SUs are allowed to access the spectrum provided that they do not interfere with the licensed users. The under-utilized spectrum bands that can be used by the SUs are called \emph{spectrum holes}~\cite{Haykin05_brain}. An ideal CR is able to efficiently detect and utilize spectrum holes.

Due to the scarce presence of active PUs and their narrow band transmission, sparse recovery methods are exploited to perform cooperative spectrum sensing \cite{Bazerque, Anese2012161group}. These approaches decompose the power spectrum density (PSD) of CRs, in terms of some appropriate bases which are related to the network topology and parameters.

\begin{figure}[b]
\vspace{-3mm}
\centering
\includegraphics[width=3.5 in,height=1.8in]{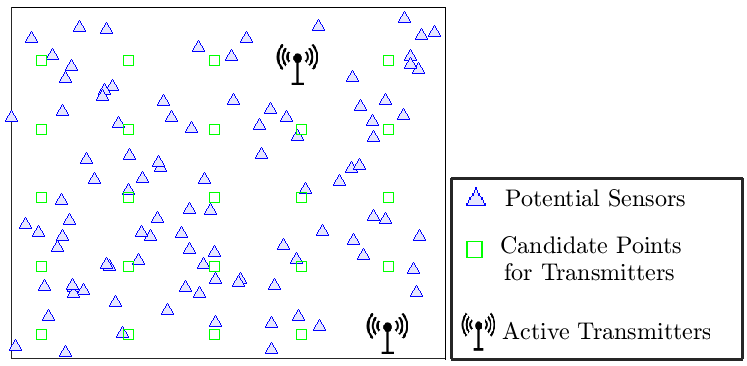}
\caption{\small{An example setup with 25 candidate points as transmitters.}}\label{setup}
\vspace{-4mm}
\end{figure}

We assume a network setup the same as that of \cite{Bazerque}. Consider $N_s$ points in a grid as the potential locations of transmitters  and  $M$ receivers in an area. The receivers receive a superposition of  transmitters signals. Figure \ref{setup} shows a setup consisting $N_s=25$ potential transmitters in which 2 of them are active and there exist $M=60$ sensors. The received signals are contaminated by channel gain and additive noise, represented by the following equation,

\begin{equation}
\label{cs_spec_sens}
\normalsize{y_m=\boldsymbol{a}_m^T \mathbf{x}+\nu_m,  \quad \quad\ \forall m=1 \ldots M},
\end{equation}
where, $\boldsymbol{a}_m$ contains the corresponding channel gains and $\nu_m$ represents noise power at the $m^{\text{th}}$ receiver. Coefficients of $\mathbf{x}$ correspond to the transmitted power at different grid points. The following problem aims to estimate $\mathbf{x}$ collaboratively using all the measurements \cite{Bazerque}.

\begin{equation}\label{reg2}
\normalsize{\hat{\mathbf{x}}=\underset{\mathbf{x}}{\text{argmin}} \|\mathbf{y}-\boldsymbol{A}\mathbf{x}\|_2^2+\gamma\|\mathbf{x}\|_1},
\end{equation}
in which, $\mathbf{y}$ and $\boldsymbol{A}$ are concatenation of ${y}_m$ and $\boldsymbol{a}_m$ respectively. Each entry of $\mathbf{x}$ determines the contribution of the $s^\text{th}$ source on the sensed data. Due to scarce presence of active transmitters and their narrow band communication, $\|\mathbf{x}\|_1$ is exploited which encourages sparsity. It should be noted that by estimating $\mathbf{x}$, we will know the location and power of transmitters at different frequency bands. These information enable us to build a radio environment map \cite{zhao2007applying}. 

\subsection{Sensor Selection Problem Statement}\label{problem_state}

Solving the sensor selection problem by evaluating the performance for each of the possible choices of $M \choose K$ is impractical unless the sizes are sufficiently small.

Suppose we want to estimate a vector $\mathbf{x}\in \mathbb{R}^N$ from $M$ linear measurements where each one is collected from a sensor, corrupted by additive noise, given by
\begin{equation}\label{sensing1}
\small{\mathbf{y}=\boldsymbol{A}\mathbf{x}+\boldsymbol{\nu}},
\end{equation}
where, $\mathbf{y}\in \mathbb{R}^M$ and $\boldsymbol{A}\in \mathbb{R}^{M\times N}$ and $\boldsymbol{\nu}$ is normally distributed with zero mean and $\sigma^2$ variance. In other words, we want to only select just $K$ rows of $\boldsymbol{A}$ to have $K$ measurements out of maximum $M$ measurements. The corresponding rows of $\boldsymbol{A}$ construct the measurement matrix, $\boldsymbol{\Phi}$, in compressed sensing literature. The maximum likelihood (ML) estimator is given by \cite{joshi_ss},
\begin{equation}\label{sol1}
\small{\hat{\mathbf{x}}_{ML}=(\boldsymbol{A}^T\boldsymbol{A})^{-1}\boldsymbol{A}^T\mathbf{y}}.
\end{equation}
The estimation error $\mathbf{x}-\hat{\mathbf{x}}$ has zero mean and the  covariance matrix is equal to
\begin{equation}\label{cov}
\small{\boldsymbol{\Sigma}_{ML}=\sigma^2(\boldsymbol{A}^T\boldsymbol{A})^{-1}}.
\end{equation}
To involve selection operator in the equations let us first write the ML solution as follows,
\begin{equation}\label{sensing}
\hat{\mathbf{x}}_{ML}={(\sum_{m=1}^M \mathbf{a}_m\mathbf{a}_m^T)}^{-1}\sum_{m=1}^M y_m\mathbf{a}_m,
\end{equation}

\normalsize{where, $\mathbf{a}_m^T$ is the $m^{\text{th}}$ row of $\boldsymbol{A}$. The estimation error is distributed in a high dimensional ellipsoid that its center is located at origin and its shape is according to the covariance matrix of error \cite{joshi_ss}. Minimization of volume of this ellipsoid (D-optimality) is the heuristic used in \cite{joshi_ss} that results in the following problem:}

\small{
\begin{equation}\label{non_convex1}
\begin{split}
&\hat{\boldsymbol{w}}=\underset{\boldsymbol{w}}{\text{argmin}} \; \text{log} \; \text{det} {(\sum_{m=1}^M w_m\mathbf{a}_m\mathbf{a}_m^T)}^{-1}, \\
&\text{subject to} \quad \|\boldsymbol{w}\|_0=K \: \text{and} \; \boldsymbol{w} \in \mathds{B}^M,
\end{split}
\end{equation}}
\normalsize{where $\boldsymbol{w}$ determines whether or not each column is involved and $\mathds{B}=\{0,1\}$.

The practical algorithms alternative to the combinatorial search are divided into two main categories, convex relaxation and greedy selection. The first approach approximates the search space to the nearest convex set and exploits convex optimization methods to solve the problem, while greedy methods gradually select suitable sensors or prune inefficient ones.}

\subsection{Matrix subset selection}
The sensor selection problem is highly related to column/row sub-matrix selection, a fundamental problem in applied mathematics. There exists many efforts in this area \cite{deshpande2010efficient,deshpande2006matrix,gu1996efficient,farahat2015greedy}. Generally, they aim at devising a computationally efficient algorithm in which the span of the selected columns/rows cover the  columns/rows space as close as possible. Mathematically, a general guarantee can be stated as one of the following forms \cite{deshpande2006matrix},
\begin{align}
\nonumber
\mathbb{E}\{ \|\boldsymbol{A}-\pi_\mathbb{T}(\boldsymbol{A})\|_F^2\}\le (K+1) \|\boldsymbol{A}-\boldsymbol{A}_K\|_F^2,\\ \nonumber
 \|\boldsymbol{A}-\pi_\mathbb{T}(\boldsymbol{A})\|_F^2\le p(K,M,N) \|\boldsymbol{A}-\boldsymbol{A}_K\|_F^2,
\end{align}

in which, $\pi_\mathbb{T}({\boldsymbol{A}})$ represents projection of rows of $\boldsymbol{A}$ on to the span of selected rows indexed by $\mathbb{T}$ set. $\mathbb{E}$ indicates expectation operator with respect to $\mathbb{T}$, i.e., all the combinatorial selection of $K$ rows of $\boldsymbol{A}$ out of $M$ are considered. Moreover, $p(K,M,N)$ is a polynomial function of the number of selected elements, the number of columns and the number of rows. $\boldsymbol{A}_K$ is the best rank-$K$ approximation of $\boldsymbol{A}$ that can be obtained by singular value decomposition. The first form suggests the distribution of potential sets for selection and it expresses an upper bound for expected value of error. The second form guarantees existence of a deterministic subset that bounds the error by a polynomial function of the parameters.     

Volume sampling is the most well-known approach to achieve the desired selection that satisfies one of the aforementioned bounds. The following theorem expresses the probabilistic form volume sampling. 

\begin{theorem}[\cite{deshpande2006matrix}]
Let $\mathbb{T}$ be a random $K-$subset of rows of a given matrix $\boldsymbol{A}$ chosen with probability 

$$
Pr(\mathbb{T})=\frac{\text{det}(\boldsymbol{A}_\mathbb{T}\boldsymbol{A}_\mathbb{T}^T)}{\sum_{|\mathbb{U}|=K} \text{det}(\boldsymbol{A}_\mathbb{U}\boldsymbol{A}_\mathbb{U}^T)}
$$

Then,
$$
\mathbb{E}\{ \|\boldsymbol{A}-\pi_\mathbb{T}(\boldsymbol{A})\|_F^2\}\le (K+1) \|\boldsymbol{A}-\boldsymbol{A}_K\|_F^2.
$$
\end{theorem}

Volume sampling considers more probability of selection for those  rows whose volume is greater.  The volume of a subset of a matrix, $\boldsymbol{A}_\mathbb{T}$, is proportional to the determinant of $\boldsymbol{A}_\mathbb{T}\boldsymbol{A}_\mathbb{T}^T$. The same heuristic criterion in (\ref{non_convex1}) aims to find the subset which has the greatest volume. It indicates the most probable subset according to volume sampling. It shows the relation of volume sampling and sensor selection in which they are solving the same problem. However this heuristic criterion is not justified for ant specific task.  

\begin{figure}[t]
\vspace{-2mm}
\centering
\includegraphics[width=2.7in]{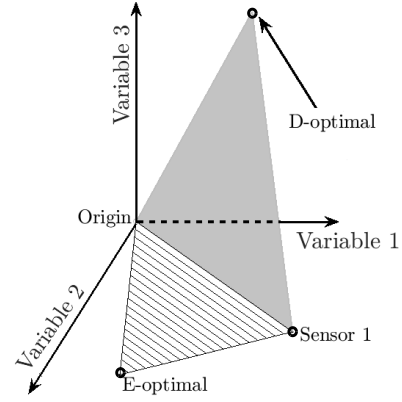}
\vspace{-5mm}
\caption{\small{Comparison of D-optimality and E-optimality for selecting 2 sensors in the 3D space. The gray area is the maximum achievable area by selecting the second sensor based on D-optimality. The shaded area is a well-shaped polygon obtained by E-optimality.}}
\label{de:opt}
\vspace{-4mm}
\end{figure}

Volume sampling and D-optimality pursue the same heuristic objective. This heuristic does not promote a well-shaped matrix for compressive sensing purposes based on RIP. However, the analysis of optimization w.r.t the RIP coefficient is not an easy task due to the columns combinatorial behavior in addition to row selection for the basic sensor selection problem. To eliminate the column combinations, we consider all of the columns and consequently we come up with an  optimization problem w.r.t the minimum eigenvalue that is known as E-optimality in the optimization literature \cite{Boyd:2004}. Assume a simple selection  from rows of $\boldsymbol{A}\in \mathbb{R}^{100\times3}$. 
Each row of $\boldsymbol{A}$, associated with a sensor, corresponds to a point in $\mathbb{R}^3$.
We are to select 2 sensors out of 100 based on D-optimality and E-optimality. Both solutions are initialized by the same sensor (sensor 1) and the criteria for the next selection varies. The D-optimal solution aims to maximize the surrounded area (gray area in Fig. \ref{de:opt}) which is vulnerable to be an ill-shaped area while, E-optimal solution comes up with a well-shaped area due to maximizing the minimum eigenvalue (shaded area in Fig. \ref{de:opt}).\footnote{The presented intuition about D-optimality and E-optimality  relates to the condition number of a matrix in linear algebra \cite{van2003numerical}. Diverged eigenvalues results in a large condition number and an ill-conditioned system of equations; accordingly, we refer to the polygon of an ill-conditioned system of equations as ill-shaped where the vertexes of shape are the rows of the matrix. On the other hand, close eigenvalues correspond to a small condition number and a well-conditioned system of equations. The corresponding polygon is referred as well-shaped in Fig \ref{de:opt}. Having well-conditioned  matrices, is a central concern in CS as evidenced  by  the  role  played  by the RIP \cite{candes2008introduction}.} 

The following  simple example illustrates the effect of E-optimality. Consider two matrices, \footnotesize{$\begin{bmatrix}

    2       & 0 \\
    0      & 0.5
\end{bmatrix}$ and $\begin{bmatrix}
  
    1       & 0 \\
    0      & 1
\end{bmatrix}$.} \normalsize{  The determinant of both matrices are equal, thus D-optimality does not favor one over the other, however, the second matrix is optimum based on E-optimality.}

As we will see in the next section, for selection of $K$ rows of $\boldsymbol{A}\in \mathbb{R}^{M\times N}$, the E-optimal criterion is equivalent to optimizing the RIP coefficient of order $N$, which is an upper bound for any arbitrary order of RIP coefficients. In the next section E-optimality will be exploited to develop a new sampling method for which its performance guarantee is analyzed. E-optimal criterion suggests optimization of an upper bound for a specific order of RIP. Moreover, in this paper we suggest a method to approximate a specific order of RIP. Based on it, a new RIP-based sensor selection algorithm is proposed.   

\normalsize{
\section{E-optimal sampling}\label{eopt_sec}
Remark 1 promotes us to develop a new matrix subset selection method that reduces the matrix to have a well-conditioned sub-matrices in the CS sense. The dominant factor of RIP constant comes from the minimum eigenvalue of the reduced matrix. It suggests to exploit the following optimization problem for sensor selection,}
\vspace{-3mm}
\small{
\begin{equation}\label{eopt_non_convex}
\begin{split}
&\hat{\boldsymbol{w}}=\underset{\boldsymbol{w}}{\text{argmin}}  \; \| {(\sum_{m=1}^M w_m\mathbf{a}_m\mathbf{a}_m^T)}^{-1}\|, \\
&\text{subject to} \quad \|\boldsymbol{w}\|_0=K \: \text{and} \; \boldsymbol{w} \in \mathds{B}^M.
\end{split}
\end{equation}}
\normalsize{
In which, $\|.\|$ denotes the spectral norm of a matrix that is defined as its maximum eigenvalue. The following lemma shows that the minimum eigenvalue is an upper bound for $\delta_S^L$.}
\begin{lemma}
For any $\boldsymbol{A}\in \mathbb{R}^{M\times N}$, the following inequality holds.
$$
1-\sigma_{\text{min}}(\boldsymbol{A}\boldsymbol{A}^T) = \delta_N^L(\boldsymbol{A}) \ge \delta_{N-1}^L(\boldsymbol{A}) \ge \cdots \ge \delta_2^L(\boldsymbol{A}).
$$
Proof: It can be concluded directly by the interlacing property of eigenvalues \cite{haemers1995interlacing}.
\label{interlace}
\end{lemma}

Lemma \ref{interlace} suggests that E-optimality, i.e., minimization of $\delta_N^L$, actually minimizes an upper bound for an arbitrary order of RIP coefficient. 

Similar to volume sampling, we design a probability of sampling according to their minimum eigenvalue. 
\begin{define}
\label{eopt_def}
Given a matrix $\boldsymbol{A}\in \mathbb{R}^{M\times N}$, \emph{E-optimal sampling} is defined as picking a subset of $\mathbb{T}$ with the following probability,
$$
Pr(\mathbb{T})=\frac{\sigma_{\text{min}}^2(\boldsymbol{A}_\mathbb{T})}{\sum_{|\mathbb{U}|=K}\sigma_{\text{min}}^2(\boldsymbol{A}_\mathbb{U})}.
$$
\end{define}
\begin{define}
\label{def_mean_rip}
Given a matrix $\boldsymbol{A}\in \mathbb{R}^{M\times N}$, $\bar{\delta}^{L}_K$ is defined as one minus the mean of minimum eigenvalues of  $\boldsymbol{A}$'s sub-matrices with $K$ columns. Mathematically, it can be expressed as follows,
$$
\bar{\delta}^{L}_K(\boldsymbol{A})=1-\mathbb{E}\{ \sigma_{\text{min}}^2(\boldsymbol{A}_\mathbb{S})\} ,
$$
in which $\mathbb{S}$ indicates a subset of $K$ columns of $\boldsymbol{A}$.
\end{define}


\begin{define}\cite{donoho2003optimally}
Given  a matrix $\boldsymbol{A}\in \mathbb{R}^{M\times N}$, the spark of $\boldsymbol{A}$ is defined as the smallest number of columns that are linearly dependent. It can be stated as follows,
\small{
$$
Spark(\boldsymbol{A})=\text{min}\; \|\textbf{x}\|_0 \; \;\; \text{s.t.} \;\;\boldsymbol{A} \textbf{x}=\textbf{0} \; \;\text{and} \; \;\textbf{x}\neq \textbf{0}. 
$$}
\end{define}
The upper bound for spark is the rank of matrix plus 1. However any linear dependencies among some columns of the matrix may decrease the spark. Based on the above definitions we present the following theorem that expresses an upper bound for projection error of E-optimal sampling. 
\begin{theorem}
Assume spark of $\boldsymbol{A}\in \mathbb{R}^{M\times N}$ is greater than $K+1$. E-optimal selection of $K$ rows implies
\small{
$$
\mathbb{E}\{\|\boldsymbol{A}-\pi_\mathbb{T}(\boldsymbol{A})\|_F^2 \}\le \frac{M-K}{C (K+1)}\frac{1-\bar{\delta}^L_{K+1}(\boldsymbol{A}^T)}{1-\bar{\delta}^L_{K}(\boldsymbol{A}^T)},
$$}
\normalsize{
where $C$ is a positive number a function of the dependencies of rows of $\boldsymbol{A}$. 
}

{proof}:
First, let us write the expansion of expected value operator according to the definition.
\small
\begin{align}
&\mathbb{E}\{\|\boldsymbol{A}-\pi_\mathbb{T}(\boldsymbol{A})\|_F^2 \} \\
&=\frac{1}{\sum_{|\mathbb{T}|=K}\sigma_{\text{min}}^2(\boldsymbol{A}_{\mathbb{T}})}\sum_{|\mathbb{T}|=K}\sigma_{\text{min}}^2(\boldsymbol{A}_{\mathbb{T}})\|\boldsymbol{A}-\pi_\mathbb{T}({\boldsymbol{A}})\|_F^2. \nonumber
\end{align}
\normalsize{
Based on the assumption on the spark of $\boldsymbol{A}$, there exist a positive constant, $\alpha$, that satisfies the following equation for every $|\mathbb{T}|=K$ and $|\mathbb{S}|=K+1$ in which $\mathbb{T}\subset \mathbb{S}$. 
}
$$
\sigma_{\text{min}}^2(\boldsymbol{A}_\mathbb{S})=\alpha \;\sigma_{\text{min}}^2(\boldsymbol{A}_\mathbb{T}) d(\boldsymbol{a}_m,\pi_{\mathbb{T}}(\boldsymbol{A})).
$$

\normalsize{
Where $\boldsymbol{a}_m$ is the innovation of $\mathbb{S}$ w.r.t $\mathbb{T}$ and $d(.,.)$ represents the Ecludian distance of a vector from a subspace. Let us take summation on all of the possible combinations,}

\small{
\begin{align}
\sum_{|\mathbb{S}|=K+1}\sigma_{\text{min}}^2 (\boldsymbol{A}_\mathbb{S})&=\sum_{|\mathbb{T}|=K}\sigma_{\text{min}}^2 (\boldsymbol{A}_\mathbb{T})\sum_{m=1}^M \alpha_m d(\boldsymbol{a}_m,\pi_{\mathbb{T}}(\boldsymbol{A}))\\& \ge C\sum_{|\mathbb{T}|=K}\sigma_{\text{min}}^2(\boldsymbol{A}_\mathbb{T})\|\boldsymbol{A}- \pi_{\mathbb{T}}(\boldsymbol{A})\|_F^2,\nonumber
\end{align}
}

\normalsize{
where $C$ is the minimum value of $\alpha_m$'s for all of the possible combinations. The assumption on the spark guarantees the existence of a positive constant. Note that the summation of the distances for all of the rows $\boldsymbol{A}$' can be stated as the Frobenius norm. Let us re-write the obtained inequality as follows,}
\vspace{-1mm}
\small{
$$
\sum_{|\mathbb{T}|=K}\sigma_{\text{min}}^2(\boldsymbol{A}_\mathbb{T})\|\boldsymbol{A}- \pi_{\mathbb{T}}(\boldsymbol{A})\|_F^2\le \frac{1}{C} \sum_{|\mathbb{S}|=K+1}\sigma_{\text{min}}^2 (\boldsymbol{A}_\mathbb{S})
$$}
Dividing  both sides of the inequality by $\sum_{|\mathbb{T}|=K}\sigma_{\text{min}}^2(\boldsymbol{A}_{\mathbb{T}})$ results in $\mathbb{E}\{ \|\boldsymbol{A}- \pi_{\mathbb{T}}(\boldsymbol{A})\|_F^2\} $ in the left side. After a simple simplification in terms of the coefficients introduced in Definition \ref{def_rip} and Definition \ref{def_mean_rip}, the right side turns into the desired expression. Please note the defined CS coefficients in (\ref{asymm_RIP}) work on the columns while we are to select some rows. $\blacksquare$
\end{theorem}

\normalsize{
E-optimal sampling implies an upper bound for the expectation of projection error in a probabilistic manner. However, we need to select some sensors deterministically. To this aim, we propose the following problem.
$$
\mathbb{S}=\underset{\mathbb{S}}{\text{argmax}} \;\lambda_{\text{min}}(A_{\mathbb{S}}A_{\mathbb{S}}^T).
$$
Algorithm  \ref{eopt_sel} shows an iterative greedy method to solve this problem. Actually, this algorithm is an approximation for the maximum likelihood estimator in which the likelihood comes from the suggested probability in Definition \ref{eopt_def}.
}

\small{
\begin{algorithm}
\caption{Greedy E-Optimal Sensor Selection}\label{eopt_sel}
\algsetup{
linenosize=\small,
linenodelimiter=:
}
\begin{algorithmic}[1]
\REQUIRE $\boldsymbol{A}$ and $K$\\
\hspace{-6mm}\textbf{Output}: The selected set $\mathbb{S}$.
\STATE \textbf{Initialization:} $\mathbb{S}$ with a random sensor
\STATE $ \text{for}\; k=1, \cdots ,K $
\STATE $\qquad \text{for}\; m=1, \cdots ,M $
\STATE $\quad\quad\qquad \mathbb{T}=\mathbb{S}\bigcup \{m\}$ 
\STATE $\quad\quad\qquad p(m)=\sigma_\text{min}(\boldsymbol{A}_\mathbb{T})$
\STATE $\qquad$end
\STATE  $s_k = \underset{m}{\text{argmax}}\;p(m)$
\STATE $\mathbb{S}=\mathbb{S}\bigcup s_k$
\STATE end
\end{algorithmic}
\end{algorithm}

\begin{table}
\centering
\caption{ \small{Complexity of different selection strategies.}}
\label{tbl_complexity}
    \begin{tabular}{  |l | l |}
    \hline
       Algorithm & Complexity \\ \hline\hline
    Convex Optimization \cite{joshi_ss} & $O(M^3)$\\ \hline
    Volume sampling \cite{deshpande2010efficient} & $O(KNM^2\text{log}M)$\\ \hline
    Greedy Submodular Selection \cite{sub_modularity}& $O(MK^3)$\\ \hline
    Greedy E-optimal selection & $O(MNK^2)$  \\ \hline
    \end{tabular}
    \vspace{-5mm}
\end{table}
\vspace{-1mm}
\normalsize

Table \ref{tbl_complexity} compares computational burden of three well-known selection methods with the proposed method. Convex relaxation is not able to work effectively for big data sets since the complexity of the algorithm grows with $M^3$ \cite{joshi_ss}. Complexity of volume sampling also depends on $M^2$. Likewise, complexity of greedy algorithms which process data one-by-one increase linearly w.r.t size of data. 

\section{Reliability Estimation and Dynamic Sensor Selection}\label{rel_section}
Collaborative sensor networks may collect redundant information which results in a larger number of sensor nodes than is needed. While, pruning unnecessary data is essential, Algorithm \ref{eopt_sel} is measurement-independent and it reduces the underlying equations of the network to shrink the equations to a well-conditioned set of sub-equations regardless of dynamic of the network. This measurement-independent approach is optimal in an averaged sense, i.e., for different possible measurements. It is appropriate for a static regime or initialization of a dynamic sensor selection. This section proposes a dynamic sensor selection framework which considers measurements for sensor selection. First of all, let us define the dynamic sensor selection systematically as follows,

\textbf{Definition 6.} (Dynamic Sensor Selection) \cite{aggarwal2011dynamic}: For a given model $\mathcal{M}$ on the data, determine set $\mathbb{S}$ such that the estimation error of the rest of sensors, $\mathbb{S}^c$, is minimized. The estimation is obtained based on the model, $\mathcal{M}$, and observed sensors, $\mathbb{S}$.

We assume the compressed sensing model (\ref{cs_spec_sens}) for power spectrum sensing as described in Section \ref{spec_sens}. Let us denote the obtained spectrum power vector by the subset $\mathbb{S}$ of sensors at time $t$ as $\boldsymbol{x}^t_\mathbb{S}$. A proper selection of $\mathbb{S}$ enables to predicting the power spectrum throughout the network's area. 


In order to keep track of the network's dynamic, we propose to sample most of the nodes in a low rate mode; while some selected  nodes should provide us with data  sampled at a high rate  enabling estimation of a high temporal resolution power spectrum map. In this framework, there is no completely switched off sensors, but we collect data from  low-sampling rate sensors to dynamically select the sensors with high sampling rate. Therefore, we have two following types of sensors in our proposed framework,
\begin{enumerate}
\normalsize{
\item \emph{High-sampling-rate selected (active) sensors}: These are a small fraction of sensors selected by an underlying sensor selection mechanism in order to  access real-time data and generate a dynamic power spectrum map. The active sensors report their sensing at rate $f_h = 1$ sample per time block.
\item \emph{Low-sampling-rate sensors}: All sensors collect and report their data in a low-rate mode, resulting in less bandwidth and power consumption. The low-rate data enables us to validate the estimated power spectrum map. The low-sampling sensors report their sensing at rate $f_{l} = \frac{1}{n_{l}}$ sample/time. I.e., 1 sample per $n_l$ time blocks is collected. It should be mentioned that the measurements from low-sampling rate sensors will not contribute in estimating $\mathbf{x}$. They will be used to determine the reliability of estimation as we will discus below.}
\end{enumerate}
\normalsize{}
The dynamic sensor selection aims to select some sensors as the active-mode set. The rest of sensors are marked as power efficient low sampling rate sensors. If the active set is selected properly, the rest of sensors can be predicted accurately by the assumed model and the active selected sensors. The ability of sensing is assumed same for all sensors and only the sensing time is different. However, different bandwidth  for sensing can be considered in a more sophisticated framework which is out of scope of this paper. Selected sensors contain sufficient information  enabling them to predict the rest of sensors by the assumed  model on the spectrum (\ref{cs_spec_sens}). Low sampling rate data may cause obsolete information vulnerable to large deviation from the model. Moreover, changes in the dynamic of network also may cause large deviations between the model's estimation and the low sampling rate data. The following expression defines a new metric  called \emph{reliability} for sensor $m$ at time $t$.
\begin{equation}\label{rel_define}
r^{(t)}_m=\frac{\text{exp}(-\sigma (t-t_m))}{1+|y_m-E^{(t)}(m,\mathbb{S})|^2}.\;\;\forall m \in \{1,\cdots,M\}
\end{equation}
in which,
$$
E^{(t)}(m,\mathbb{S})=\boldsymbol{a}_m^T\boldsymbol{x}^{(t-1)}(\mathbb{S})
$$
In (\ref{rel_define}), $\boldsymbol{x}^{(t-1)}(\mathbb{S})$ is the estimation of power propagation at time $t-1$ based on collected data from active sensors indexed by $\mathbb{S}$. Moreover, $E^{(t)}(m,\mathbb{S})$ is the estimation of the measurement of $m^{\text{th}}$ sensor at time $t$. $\sigma$ is a temporal forgetting factor. $t_m$ is the last time that sensor $m$ is sampled and the corresponding measurement is $y_m$. The reliability of each sensor consists of two terms. The numerator indicates how fresh is our observation. Obsolete data results in unreliable observation. The denominator shows the power of model for estimation of unseen regions. Accurate estimation of the observation of sensor $m$ using the active demonstrates that the sensor $m$ has a reliable sensing. 
The proposed dynamic sensor selection framework is illustrated in Fig. \ref{dynamic_SS}. We propose to consider the reliability of sensors in the sensor selection procedure in order to determine a proper subset which is able to compensate large model's error for the low-rate sampled sensors. Mathematically speaking, the static E-optimal sensor selection algorithm is modified as follows,

\begin{figure}[t]
\centering
\vspace{-2mm}
\includegraphics[width=2.7 in,angle=0]{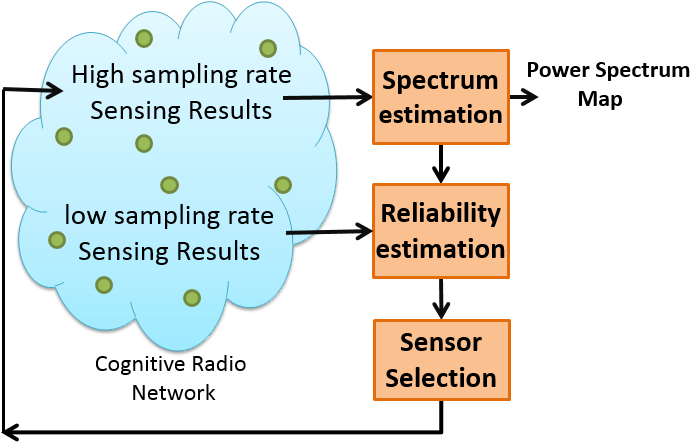}
\caption{\small{The main framework of the proposed reliability based sensor selection.}}
\vspace{-5mm}
\label{dynamic_SS}
\end{figure}

\begin{equation}
\label{DSS}
\mathbb{S}=\underset{|\mathbb{S}|\le K}{\text{argmax}} \;\lambda_{\text{min}}(A_{\mathbb{S}}A_{\mathbb{S}}^T)+\gamma \|u_\mathbb{S}\|_2^2,
\end{equation}
in which, $\gamma $ is the regularization parameter and $u_m=r^{-1}_m$ represents unreliability and $u_\mathbb{S}$ is the sub-vector of $\boldsymbol{u}$ indexed by set $\mathbb{S}$. The superscript $(t)$ is removed due to simplicity of notation. It means we are looking for unreliable sensors to select them for the next time slot in order to compensate the model's error.  

\section{Optimization and Complexity}\label{opt_sec}

In order to cast the dynamic sensor selection (\ref{DSS}) in a tractable formulation, first let us rewrite the minimum eigenvalue as the following problem.

\small{
\begin{equation}\label{RIP_problem}
\lambda_{\text{min}}(\boldsymbol{A})=\text{min} \|\boldsymbol{A}\mathbf{x}\|_2^2 \quad \text{s.t.} \; \|\mathbf{x}\|_2=1.
\end{equation}}
\normalsize{Problem (\ref{DSS}) can be written in the following form,}
\small{
\begin{align}\label{DSS2}
\boldsymbol{W}^{(t)}&=\underset{W}{\text{argmax}} \;\underset{x}{\text{min}} \|\boldsymbol{WA}\mathbf{x}\|_2^2 +\gamma\|\boldsymbol{Wu}^{(t)}\|_2^2\quad \text{s.t.}\\  & \|\mathbf{x}\|_2=1 \;, W_{ij}\in\{0,1\}, \|\mathbf{w}_k\|_0=1 \text{ and }\|\mathbf{w}^m\|_0\le 1. \nonumber
\end{align}}
\vspace{-3mm}

\normalsize{In which $\boldsymbol{W}\in \mathbb{R}^{K\times M}$ reduces the matrix $\boldsymbol{A}\in \mathbb{R}^{M\times N}$ by some selected rows. $\mathbf{w}_k$ represents the $k^{\text{th}}$ row of $\boldsymbol{W}$ and $\mathbf{w}^m$} indicates the $m^{\text{th}}$ column of $\boldsymbol{W}$. \normalsize{ The last constraint $\small{\|\mathbf{w}^m\|_0\le 1}$ avoids repetitive selection of the same row (sensor). This problem implies eigenvalue optimization over combination of rows of $\boldsymbol{A}$ that it is shown to be NP-hard \cite{Tillmann:2014:CCR:2689743.2690742}. Accordingly, we propose a greedy algorithm to solve (\ref{DSS2}).
}

\normalsize{
Algorithm \ref{sens_sel1} shows the steps of our proposed greedy algorithm to solve the obtained optimization problem. This algorithm optimizes the reduction matrix row-by-row where the reliability of the non-selected sensors are being considered. Assume the algorithm aims to select a new sensor at the $k^{\text{th}}$ iteration. Up to current iteration, $k-1$ sensors already are selected. The algorithm evaluate the non-selected sensors one-by-one in order to find the sensor that maximizes the objective function. The objective function is a weighted summation of the minimum eigenvalue of the restricted set of rows (sensors) and their corresponding unreliability weights. To evaluate each sensor we need to compute the most dominant $k$ eigen components which implies performing singular value decomposition (SVD). However, truncated SVD up to the $k^{\text{th}}$ component will be sufficient. An online algorithm is proposed that observes the non-selected sensors with a low sampling rate as depicted in Fig. \ref{dynamic_SS}. In each sequence, the observed set of sensors is updated as well as their corresponding reliability weights. The first step to update the reliability is estimating the propagation using only the current active sensors. To this aim the following problem must be solved.  
\begin{equation}
\label{x_lasso}
\boldsymbol{x}^{(t)}(\mathbb{S})=\underset{\boldsymbol{x}}{\text{argmin}}\; \|\boldsymbol{W}^{(t)}(\boldsymbol{y}-A\boldsymbol{x})\|_2^2 + \lambda_{LASSO} \|\boldsymbol{x}\|_1.
\end{equation}
\vspace{-3mm}
\small{
\begin{algorithm}[t]
\caption{\small{Reliable E-optimal Sensor Selection}}\label{sens_sel1}
\algsetup{
linenosize=\small,
linenodelimiter=:
}
\begin{algorithmic}
\vspace{-1mm}
\REQUIRE $\boldsymbol{A}$, $S$, $K$ and $\boldsymbol{r}$\\
\hspace{-3mm}\textbf{Output}: The selected set $\mathbb{S}$ and reduction matrix $\boldsymbol{W}$.
\STATE $ \textbf{Initialization:  } \boldsymbol{W}=\boldsymbol{0} \in \mathbb{R}^{K\times M}$ and $\mathbb{S}=\emptyset$
\STATE $ \text{for}\; k=1, \cdots ,K $ (Optimization of the $k^{\text{th}}$ row of $\boldsymbol{W}$)
\STATE $\quad\qquad \text{for}\; \forall m\in \mathbb{S}^c$
\STATE $\quad\quad\qquad$SVD:  $\boldsymbol{A}(\mathbb{S}\bigcup m,:)=\boldsymbol{V}^T\boldsymbol{\Lambda U}$
\STATE $\quad\quad\qquad\mathbf{x}^*= U(:,k)$
\STATE $\quad\quad\qquad p(m)=\|\boldsymbol{A}\mathbf{x}^*\|_2^2+\gamma u(m)$
\STATE $\quad\qquad$end
\STATE  $\quad s_k = \underset{m}{\text{argmax}}\;p(m)$
\STATE $\quad \mathbb{S}=\mathbb{S}\bigcup s_k$ and $\boldsymbol{W}_{k,s_k}=1$
\STATE end for

\end{algorithmic}
\vspace{-1mm}
\end{algorithm}

\normalsize{Here $\lambda_{LASSO}$ regularizes  sparsity and $\boldsymbol{W}$ indicates the reduction matrix to the selected set  $\mathbb{S}$. Those sensors whose measurements are matched with the estimated power density map are marked as reliable. A consistent definition is proposed in (\ref{rel_define}) which considers the deviation of actual measurements from the estimation of the model as a metric for reliability. The subscript $t$ is removed in Algorithm 2 for simplification. Algorithm \ref{sens_sel2} shows the overall process of spectrum sensing using the selected sensors.

\normalsize{}
The bottleneck of complexity order of Algorithm \ref{sens_sel1} at the $k^{\text{th}}$ iteration is performing a truncated singular value decomposition to obtain the first $k$ eigen components. Thus, the complexity of the algorithm in the $k^{\text{th}}$ iteration will be $O(kMN^2)$ \cite{holmes2007fast}. Therefore, selection of $K$ sensors implies complexity order of $O(K^2MN^2)$.   

\small{
\begin{algorithm}
\caption{\small{Spectrum Sensing using Dynamic Sensor Selection}}\label{sens_sel2}
\algsetup{
linenosize=\small,
linenodelimiter=:
}
\begin{algorithmic}[2]
\vspace{-1mm}
\REQUIRE $\boldsymbol{A}$, $S$, $K$, $\lambda$, $f_l$ and $\lambda_{LASSO}$. \\
\hspace{-6mm}\textbf{Output}: Power spectrum for each time $\boldsymbol{x}^{(t)}$.
\STATE $ \textbf{Initialization: }   \mathbb{S}=$  Output of Algorithm 1 and $\boldsymbol{x}(\mathbb{S})=$ Result of Problem (\ref{x_lasso})
\STATE  for a new time block ($t$)
\STATE $\quad$ sample $M\times f_l$ sensors
\STATE  $\quad$ Update $t_m=t$ for the sensed sensors
\STATE  $\quad$ Update reliability using (\ref{rel_define})  
\STATE $\quad $ $\mathbb{S}^{(t)}=$ Output of Algorithm \ref{sens_sel1}
\STATE $\quad$  $\boldsymbol{x}^{(t)}(\mathbb{S}^{(t)})=$ Result of Problem (\ref{x_lasso})
\STATE end for
\end{algorithmic}
\end{algorithm}
\vspace{-1mm}
}
\normalsize
\vspace{-1mm}
\normalsize{
\section{Experimental Results}\label{Experimental}

The simulations are performed for collaborative spectrum sensing. The setup for generating data are employed from \cite{Bazerque}. Our goal is to estimate vector $\mathbf{x}$ that indicates transmitted spectrum power at some candidate points. 

For the first simulation suppose we have potentially $300$ sensors and they are estimating an $\textbf{x} \in \mathbb{R}^{36}$ that has only $5$ active transmitters. The location of sensors are derived from a uniform distribution and the active transmitters are selected randomly and the results are averaged for 200 different realizations. The following linear measurements are sensed by sensors $m=1 \ldots M$,
$$
\normalsize{y_m=\boldsymbol{a}_m^T \mathbf{x}+\nu_m },
$$
where, $\nu_m$ indicates additive white Gaussian noise.
$a_{ms}$ shows the $s^{\text{th}}$ entry of $\boldsymbol{a}_m$ is the channel gain between the $m^{\text{th}}$ sensor and the $s^{\text{th}}$ potential source. The channel gain between two points is assumed by one over squared distance of two points. Since, the ability of sensors is considered the same over spectrum, thus the simulations are performed for a single spectrum band. The same procedure can be performed for multi-band spectrum regime independently. Figure \ref{exprmn1} shows the performance of different static algorithms versus the number of selected sensors. Static refers to measurement-independent methods. In this experiment the SNR is set to +20dB. Successful recovery is defined as true estimation of the support of sparse vector using the measurements. Problem (\ref{x_lasso}) is solved $200$ for each algorithm. The Sparse solution is obtained using the iterative re-weighted least square algorithm \cite{IRLS}. As it can be seen in Fig. \ref{exprmn1}, E-optimal based sensor selection has the best performance.

\begin{figure}[b]
\centering
\vspace{-5mm}
\includegraphics[width=3 in,height=1.8 in]{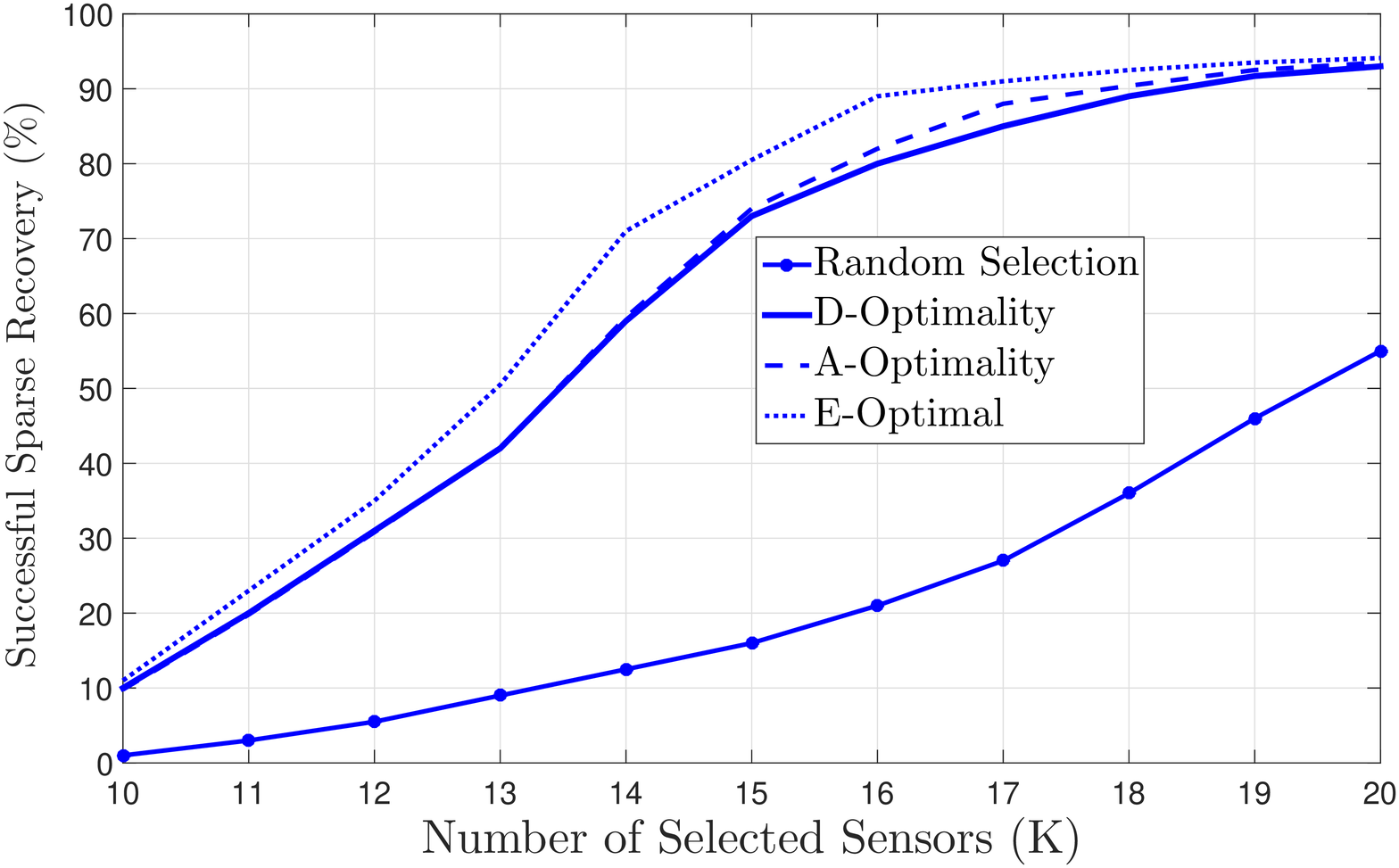}
\caption{\small{Performance of different static sensor selection algorithms in terms of number of selected sensors.}}\label{exprmn1}
\vspace{-6mm}
\end{figure}

\begin{figure}[t]
\centering
\includegraphics[width=3 in,angle=0]{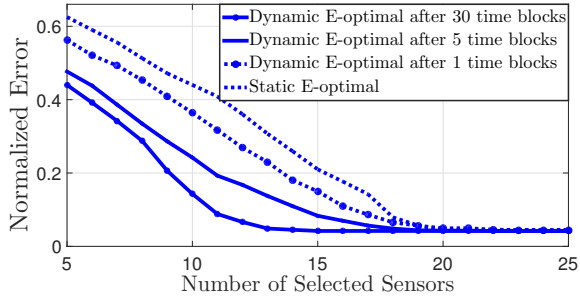}
\caption{\small{Performance of static and dynamic E-optimal-based sensor selection algorithms vs. the number of selected sensors.}}
\label{exprmn_mse}
\vspace{-5mm}
\end{figure}


Fig. \ref{exprmn_mse} exhibits the effect of involving  reliability on the static sensor selection.  Suppose there are $300$ potential sensors and the low-sampling rate is set equal to $\frac{1}{30}$. It means in each time block $10$ new measurements contribute to construct the reliability weights (\ref{rel_define}). Observation of new measurements of one time block makes an improvement in normalized estimation error; similarly, usage of $5$ time blocks significantly improves the performance to be close to the estimation after $30$ time blocks in which all the sensors are observed. The forgetting factor is set to $0$ as the state of network is not changed during observation of $30$ time blocks. Thus, aggregating the measurements without the forgetting factor is optimum. The normalized error,
$
\small{\|\textbf{x}^*-{\textbf{x}(\mathbb{S})}\|_2/{\|\textbf{x}^*\|_2}}
$, is defined as the criterion for performance, where, $\textbf{x}^*$ is the ground truth solution. 

Fig. \ref{sensors_compare} visualizes the error of spectrum sensing in the area of network for the setup of Fig. \ref{exprmn_mse}. We are to choose 8 sensors.

Fig. 7 shows that the error of estimation is significantly decreased by setting $\gamma=0.7$ for the setup of Fig. \ref{exprmn_mse}. However, an efficient value of $\gamma$ depends on the problem setup and should be tuned.  Setting $\gamma=0$ is equivalent to the static E-optimal sensor selection. Simulation shows the proposed reliable sensor selection performs better than the static sensor selection for a relatively wide range of $\gamma$, i.e., the problem is not very sensitive to well-tuning of this parameter.
\begin{figure}[t]
\centering
\includegraphics[width=3.5 in,angle=0]{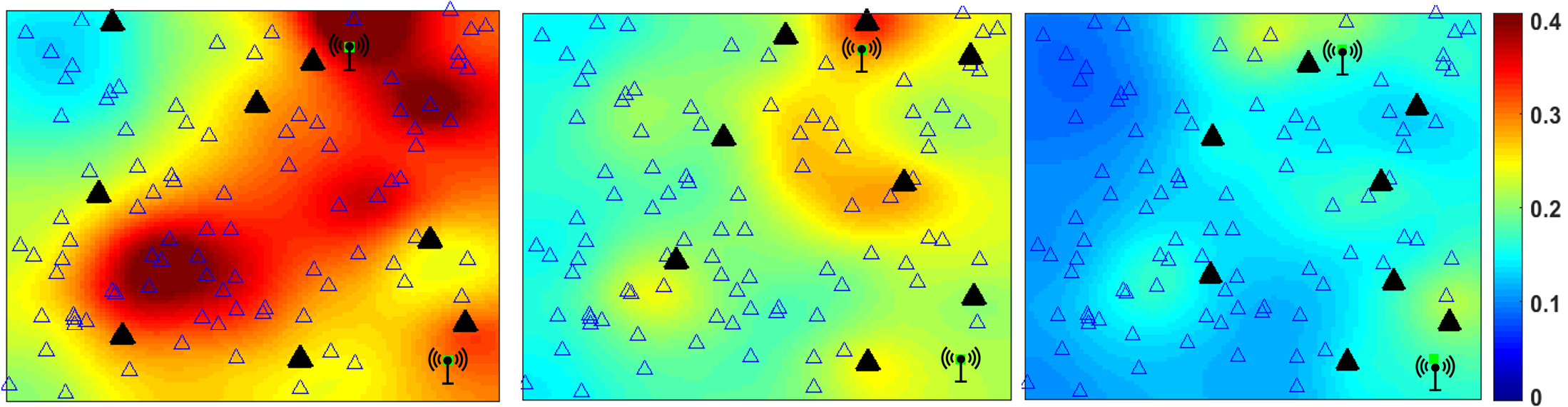}
\caption{\small{The error of estimated spectrum in the area of interest. (Left) E-optimal, Algorithm 1. (Middle) Reliable E-optimal, Algorithm 2 after sensing in one time block. (Right) Reliable E-optimal, Algorithm 2 while all the sensors are sensed after 30 time blocks. $\gamma$ is assumed $0.7$}}\label{sensors_compare}
\vspace{-4mm}
\end{figure}

\begin{figure}[b]
\centering
\vspace{-2mm}
\includegraphics[width=1.7 in,angle=0]{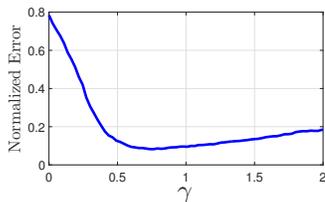}
\vspace{-3mm}
\small{\caption{MSE error versus different values of $\gamma$.}}
\vspace{-2mm}
\label{lambda}
\end{figure}


\begin{figure}
\centering
\includegraphics[width=3.5 in,angle=0]{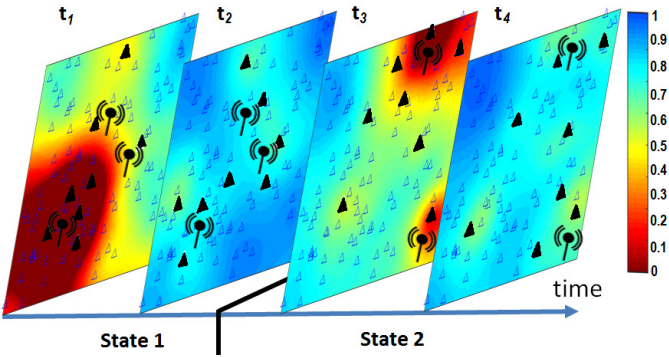}
\caption{\small{Reliability maps of $4$ time blocks illustrate how the proposed framework evolves in time in order to select adapted sensors to the dynamic of network after state transition. Sensors within unreliable (red) areas have more chance of selection. }}\label{dynamic_vis}
\vspace{-6mm}
\end{figure}
\begin{figure}[t]
\centering
\begin{subfigure}{0.5\textwidth}
\centering
\includegraphics[width=3 in,angle=0]{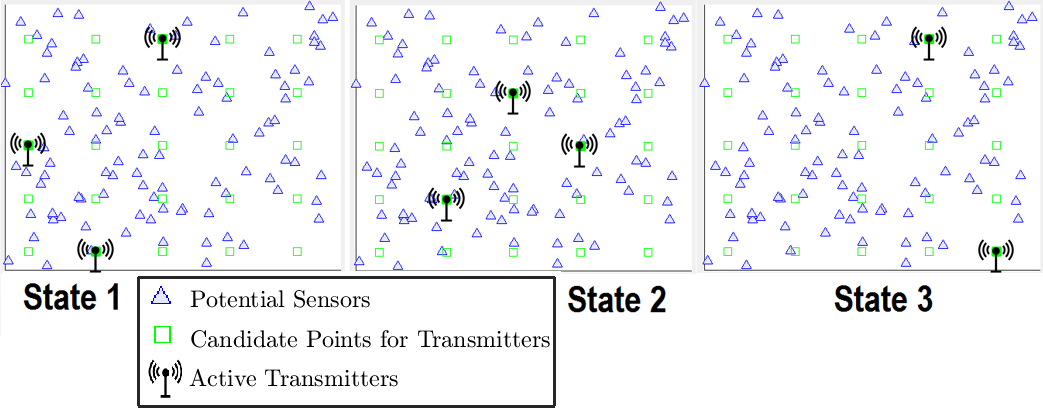}
\vspace{-2mm}
\caption{ }\label{states}
\end{subfigure}
\begin{subfigure}{0.5\textwidth}
\centering
\includegraphics[width=3.5 in,height=1.7 in]{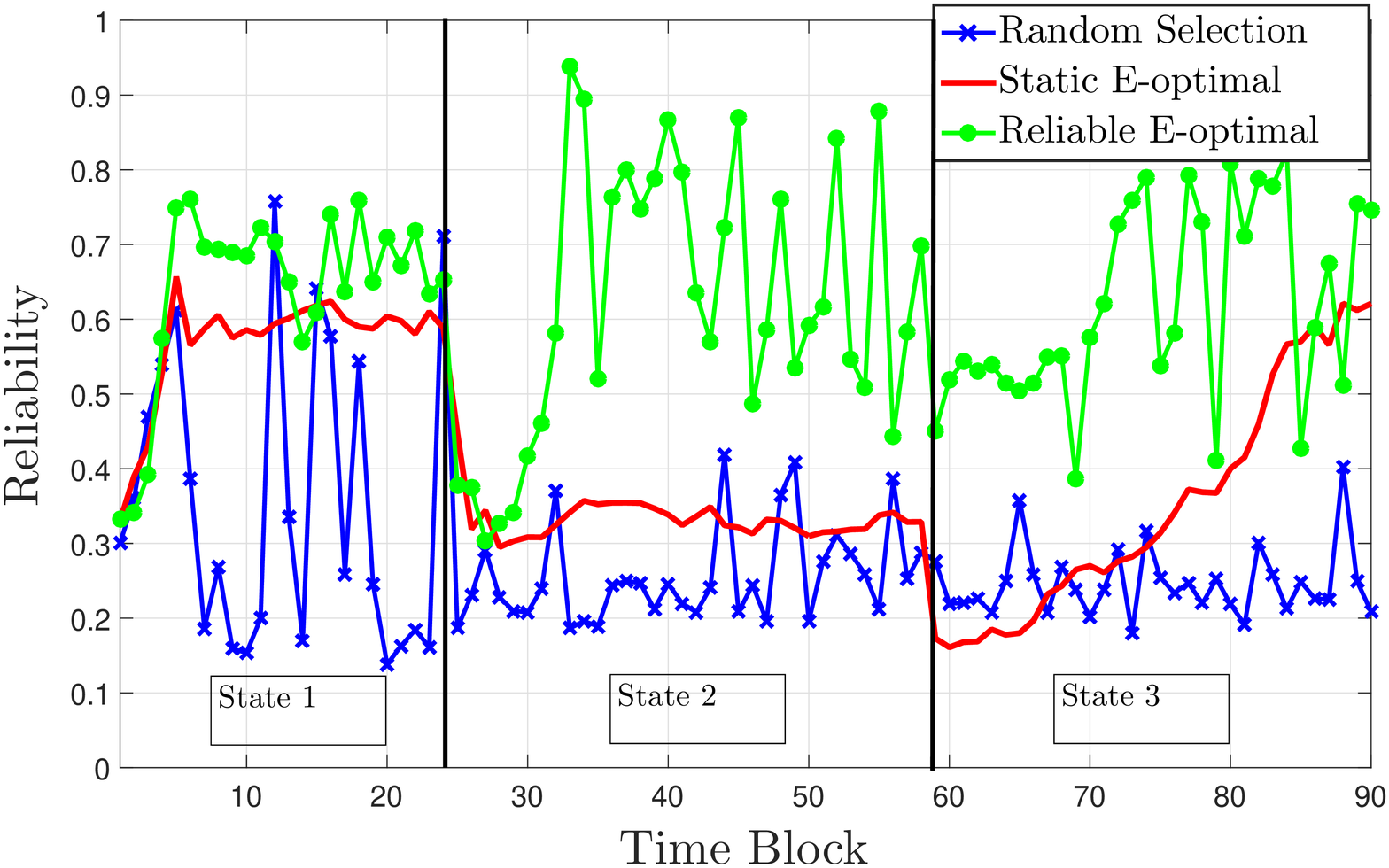}
\vspace{-2mm}
\caption{}\label{reliablity_fig}
\end{subfigure}
\vspace{-2mm}
\caption{(a)\small{ A dynamic network with $3$ states for the location of active PUs. The shaded blue squares represent  active PUs. (b) The effect of reliable sensor selection for compensation of the model error in the reliable sensor selection procedure.}}
\vspace{-5mm}
\end{figure}
\begin{figure}[t]
\centering
\includegraphics[width=3.5 in,height=1.6 in]{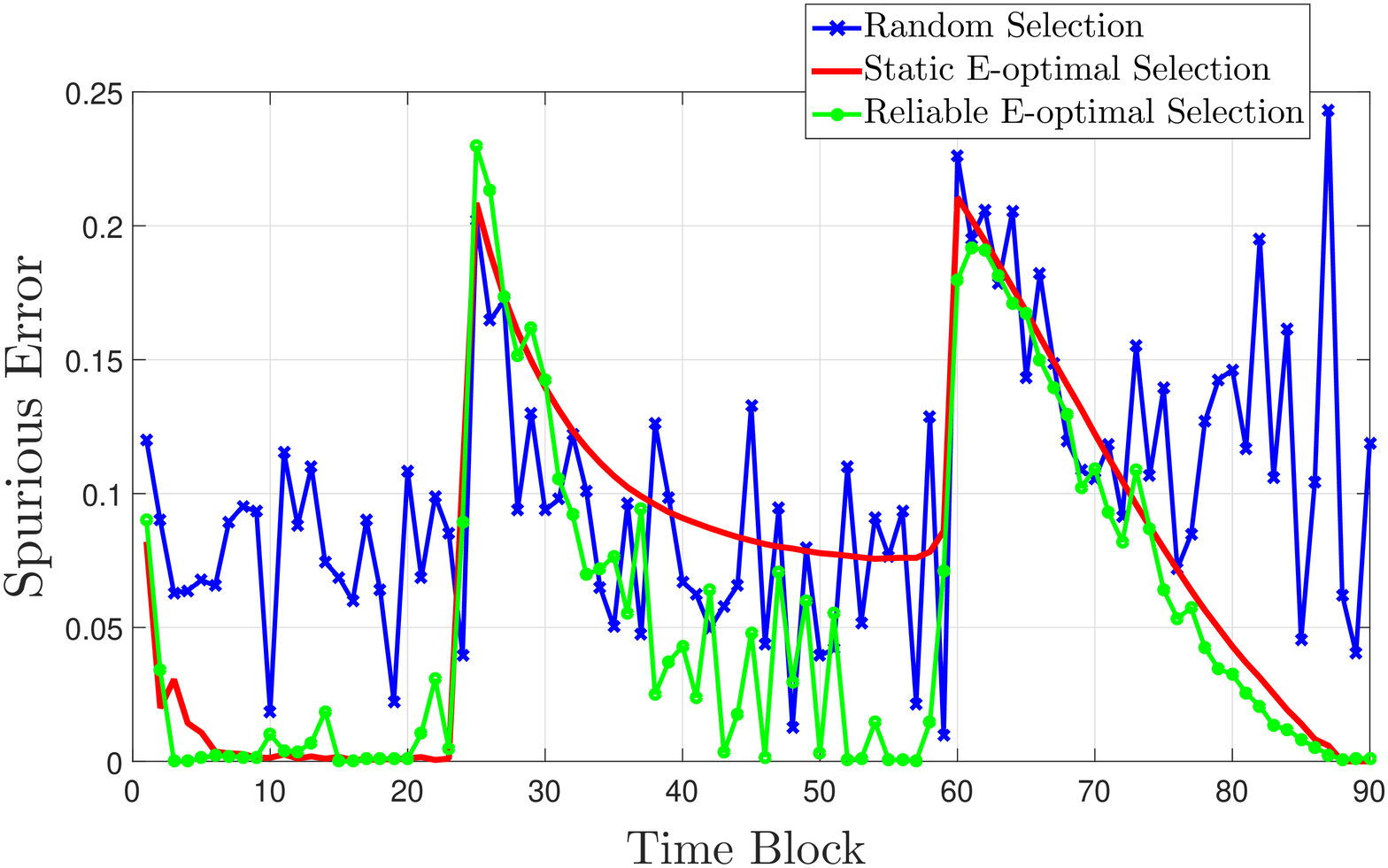}
\vspace{-3mm}
\caption{\small{The effect of reliable sensor selection for compensation of the model error in the reliable sensor selection procedure. }}\label{spurious}
\vspace{-6mm}
\end{figure}

In addition to power spectrum map, the proposed framework is able to generate a new network profile which can provides us trustworthy of the estimated spectrum for each point of the network. We call this side output \emph{reliability map}. Interpolation of the estimated reliability of sensors throughout the network's area, generates the reliability map. Fig. \ref{dynamic_vis} visualizes the temporal effect of dynamic sensor selection using the reliability map. Unreliable areas are indicated by red and blue areas represent reliable estimation of spectrum. Reliable sensor selection aims to compensate unreliability by considering more chance for red regions. In the next time slot the error for those regions are compensated. In this figure, each state of the network corresponds to a specific set of active PUs.

Fig. \ref{states} shows the location of active PUs for a dynamic network with $3$ states. There are $90$ time blocks and the state of network is changed in blocks $24$ and $59$. The forgetting factor is set to $0.1/(\Delta T)$ in which $\Delta T$ is the time difference of two consecutive time blocks. Fig. \ref{reliablity_fig} and Fig. \ref{spurious} show the performance of sensor selection in terms of average network reliability and the spurious error of spectrum sensing which is defined by
$
\|\hat{\boldsymbol{x}}_\text{spurious}\|_1=\sum_{i\not\in \boldsymbol{x}^* \;\text{support}}| {x(\mathbb{S})}_i|.
$ As it can be seen, the reliability is increased and the undesired power propagation is decreased  by exploiting the dynamic framework, especially for the second state. 
\vspace{-1mm}
\section{Conclusion}\label{conc}
\vspace{-1mm}
The problem of sensor selection is considered and its relation to existing work on matrix subset selection is elaborated. We developed a new subset selection method as an extension of the well-known volume sampling. Our criteria is based on E-optimality, which is in favor of compressive sensing theory. We extended the static E-optimal sensor selection to a dynamic sensor selection method that exploits the measurements in an online manner. The experimental results indicate the efficiency of our suggested sensor selection algorithm in cognitive radio networks' spectrum sensing.


\small{
\balance
\bibliographystyle{IEEEtran}
\bibliography{ref}

\begin{thebibliography}{10}
\providecommand{\url}[1]{#1}
\csname url@samestyle\endcsname
\providecommand{\newblock}{\relax}
\providecommand{\bibinfo}[2]{#2}
\providecommand{\BIBentrySTDinterwordspacing}{\spaceskip=0pt\relax}
\providecommand{\BIBentryALTinterwordstretchfactor}{4}
\providecommand{\BIBentryALTinterwordspacing}{\spaceskip=\fontdimen2\font plus
\BIBentryALTinterwordstretchfactor\fontdimen3\font minus
  \fontdimen4\font\relax}
\providecommand{\BIBforeignlanguage}[2]{{%
\expandafter\ifx\csname l@#1\endcsname\relax
\typeout{** WARNING: IEEEtran.bst: No hyphenation pattern has been}%
\typeout{** loaded for the language `#1'. Using the pattern for}%
\typeout{** the default language instead.}%
\else
\language=\csname l@#1\endcsname
\fi
#2}}
\providecommand{\BIBdecl}{\relax}
\BIBdecl

\bibitem{joshi_ss}
S.~Joshi and S.~Boyd, ``Sensor selection via convex optimization,''
  \emph{Signal Processing, IEEE Transactions on}, vol.~57, no.~2, pp. 451--462,
  Feb 2009.

\bibitem{sub_modularity}
M.~Shamaiah, S.~Banerjee, and H.~Vikalo, ``Greedy sensor selection: Leveraging
  submodularity,'' in \emph{Decision and Control (CDC), 2010 49th IEEE Conf.
  on}, Dec 2010, pp. 2572--2577.

\bibitem{shirazi}
M.~Shirazi, A.~Sani, and A.~Vosoughi, ``Sensor selection and power allocation
  via maximizing bayesian fisher information for distributed vector
  estimation,'' in \emph{2017 51th Asilomar Conference on Signals, Systems and
  Computers}, Nov 2017.

\bibitem{Boyd:2004}
S.~Boyd and L.~Vandenberghe, \emph{Convex Optimization}.\hskip 1em plus 0.5em
  minus 0.4em\relax New York, NY, USA: Cambridge University Press, 2004.

\bibitem{cs_book}
S.~Foucart and H.~Rauhut, \emph{A Mathematical Introduction to Compressive
  Sensing}.\hskip 1em plus 0.5em minus 0.4em\relax Birkh\"{a}user, 2013.

\bibitem{cs_app_iamge}
B.~Amizic, L.~Spinoulas, R.~Molina, and A.~Katsaggelos, ``Compressive blind
  image deconvolution,'' \emph{Image Processing, IEEE Trans. on}, vol.~22,
  no.~10, pp. 3994--4006, Oct 2013.

\bibitem{cs_app_comm}
X.~Guan, Y.~Gao, J.~Chang, and Z.~Zhang, ``Advances in theory of compressive
  sensing and applications in communication,'' in \emph{Instrumentation,
  Measurement, Computer, Communication and Control, 2011 First International
  Conference on}, Oct 2011.

\bibitem{cs_app_ss}
Y.~Gwon, H.~Kung, and D.~Vlah, ``Compressive sensing with optimal sparsifying
  basis and applications in spectrum sensing,'' in \emph{Global Communications
  Conference (GLOBECOM), 2012 IEEE}, Dec 2012, pp. 5386--5391.

\bibitem{omp}
T.~Cai and L.~Wang, ``Orthogonal matching pursuit for sparse signal recovery
  with noise,'' \emph{Information Theory, IEEE Transactions on}, vol.~57,
  no.~7, pp. 4680--4688, July 2011.

\bibitem{slzero}
H.~Mohimani, M.~Babaie-Zadeh, and C.~Jutten, ``A fast approach for overcomplete
  sparse decomposition based on smoothed $\ell_0$ norm,'' \emph{Signal
  Processing, IEEE Transactions on}, vol.~57, no.~1, pp. 289--301, Jan 2009.

\bibitem{basispursut}
S.~S. Chen, D.~L. Donoho, and M.~A. Saunders, ``Atomic decomposition by basis
  pursuit,'' \emph{SIAM JOURNAL ON SCIENTIFIC COMPUTING}, vol.~20, pp. 33--61,
  1998.

\bibitem{dic_learning}
Q.~Geng and J.~Wright, ``On the local correctness of $\ell_1$ minimization for
  dictionary learning,'' in \emph{Information Theory (ISIT), 2014 IEEE
  Symposium on}, June 2014, pp. 3180--3184.

\bibitem{Donoho01uncertainty}
D.~L. Donoho and X.~Huo, ``Uncertainty principles and ideal atomic
  decomposition,'' \emph{IEEE Transactions on Information Theory}, 2001.

\bibitem{Cohen09compressedsensing}
A.~Cohen, W.~Dahmen, and R.~Devore, ``Compressed sensing and best k-term
  approximation,'' \emph{J. Amer. Math. Soc}, 2009.

\bibitem{donoho2003optimally}
D.~L. Donoho and M.~Elad, ``Optimally sparse representation in general
  (nonorthogonal) dictionaries via $\ell_1$ minimization,'' \emph{Proceedings
  of the National Academy of Sciences}, vol. 100, no.~5, pp. 2197--2202, 2003.

\bibitem{Candes_RIP}
E.~Candes and T.~Tao, ``Decoding by linear programming,'' \emph{Information
  Theory, IEEE Transactions on}, vol.~51, no.~12, pp. 4203--4215, Dec 2005.

\bibitem{Candes_RIP2}
E.~Candes, ``The restricted isometry property and its implications for
  compressed sensing,'' \emph{C. R. Academie des Sciences}, no. 356, pp.
  689--592, 2008.

\bibitem{Davies_RIP3}
M.~Davies and R.~Gribonval, ``Restricted isometry constants where $\ell_p$
  sparse recovery can fail for $0< p \leq 1$,'' \emph{Information Theory, IEEE
  Transactions on}, vol.~55, no.~5, May 2009.

\bibitem{Blanchard:2011}
J.~D. Blanchard, C.~Cartis, and J.~Tanner, ``Compressed sensing: How sharp is
  the restricted isometry property?'' \emph{SIAM Rev.}, vol.~53, no.~1, pp.
  105--125, Feb. 2011.

\bibitem{Akyildiz06_xG}
I.~F. Akyildiz, W.-Y. Lee, M.~C. Vuran, and S.~Mohanty, ``Next
  generation/dynamic spectrum access/cognitive radio wireless networks: A
  survey,'' \emph{Computer Networks}, vol.~50, no.~13, pp. 2127 -- 2159, May
  2006.

\bibitem{Dynamic_Spectrum_Access4205091}
Q.~Zhao and B.~Sadler, ``A survey of dynamic spectrum access,'' \emph{Signal
  Processing Magazine, IEEE}, vol.~24, no.~3, pp. 79--89, May 2007.

\bibitem{Haykin05_brain}
S.~Haykin, ``Cognitive radio: brain-empowered wireless communications,''
  \emph{IEEE Journal on Selected Areas in Communications}, vol.~23, no.~2, pp.
  201--220, Feb. 2005.

\bibitem{Bazerque}
J.~Bazerque and G.~Giannakis, ``Distributed spectrum sensing for cognitive
  radio networks by exploiting sparsity,'' \emph{Signal Processing, IEEE
  Transactions on}, vol.~58, no.~3, pp. 1847--1862, March 2010.

\bibitem{Anese2012161group}
E.~Dall'Anese, J.~A. Bazerque, and G.~B. Giannakis, ``Group sparse lasso for
  cognitive network sensing robust to model uncertainties and outliers,''
  \emph{Physical Communication}, vol.~5, no.~2, pp. 161 -- 172, 2012,
  compressive Sensing in Communications.

\bibitem{zhao2007applying}
Y.~Zhao, L.~Morales, J.~Gaeddert, K.~Bae, J.-S. Um, and J.~Reed, ``Applying
  radio environment maps to cognitive wireless regional area networks,'' in
  \emph{New Frontiers in Dynamic Spectrum Access Networks, 2007. DySPAN 2007.
  2nd IEEE International Symposium on}, April 2007, pp. 115--118.

\bibitem{deshpande2010efficient}
A.~Deshpande and L.~Rademacher, ``Efficient volume sampling for row/column
  subset selection,'' in \emph{Foundations of Computer Science (FOCS), 2010
  51st Annual IEEE Symposium on}.\hskip 1em plus 0.5em minus 0.4em\relax IEEE,
  2010, pp. 329--338.

\bibitem{deshpande2006matrix}
A.~Deshpande, L.~Rademacher, S.~Vempala, and G.~Wang, ``Matrix approximation
  and projective clustering via volume sampling,'' in \emph{Proceedings of the
  seventeenth annual ACM-SIAM symposium on Discrete algorithm}, 2006, pp.
  1117--1126.

\bibitem{gu1996efficient}
M.~Gu and S.~C. Eisenstat, ``Efficient algorithms for computing a strong
  rank-revealing qr factorization,'' \emph{SIAM Journal on Scientific
  Computing}, vol.~17, no.~4, pp. 848--869, 1996.

\bibitem{farahat2015greedy}
A.~Farahat, A.~Elgohary, and A.~Ghodsi, ``Greedy column subset selection for
  large-scale data sets,'' \emph{Knowledge and Information Systems}, vol.~45,
  no.~1, pp. 1--34, 2015.

\bibitem{van2003numerical}
P.~Van~Dooren, ``Numerical linear algebra for signal, systems and control,''
  \emph{Draft notes prepared for the Graduate School in Systems and Control},
  vol. 250, 2003.

\bibitem{candes2008introduction}
E.~J. Cand{\`e}s and M.~B. Wakin, ``An introduction to compressive sampling,''
  \emph{IEEE signal processing magazine}, vol.~25, no.~2, pp. 21--30, 2008.

\bibitem{haemers1995interlacing}
W.~H. Haemers, ``Interlacing eigenvalues and graphs,'' \emph{Linear Algebra and
  its applications}, vol. 226, pp. 593--616, 1995.

\bibitem{aggarwal2011dynamic}
C.~C. Aggarwal, Y.~Xie, and P.~S. Yu, ``On dynamic data-driven selection of
  sensor streams,'' in \emph{Proceedings of the 17th ACM SIGKDD international
  conference on Knowledge discovery and data mining}.\hskip 1em plus 0.5em
  minus 0.4em\relax ACM, 2011, pp. 1226--1234.

\bibitem{Tillmann:2014:CCR:2689743.2690742}
A.~M. Tillmann and M.~E. Pfetsch, ``The computational complexity of the
  restricted isometry property, the nullspace property, and related concepts in
  compressed sensing,'' \emph{IEEE Trans. Inf. Theor.}, vol.~60, no.~2, pp.
  1248--1259, Feb. 2014.

\bibitem{holmes2007fast}
M.~Holmes, A.~Gray, and C.~Isbell, ``Fast svd for large-scale matrices,'' in
  \emph{Workshop on Efficient Machine Learning at NIPS}, vol.~58, 2007, pp.
  249--252.

\bibitem{IRLS}
R.~Chartrand and W.~Yin, ``Iteratively reweighted algorithms for compressive
  sensing,'' in \emph{Acoustics, Speech and Signal Processing, 2008. ICASSP
  2008. IEEE International Conference on}, March 2008, pp. 3869--3872.

\end{thebibliography}
}
\end{document}